\documentclass[12pt,bezier]{amsart}
\usepackage{amsfonts,latexsym,amscd}

\newtheorem{theorem}{Theorem}[section]

\newtheorem{proposition}[theorem]{Proposition}
\newtheorem{corollary}[theorem]{Corollary}
\newtheorem{lemma}[theorem]{Lemma}

\newcommand{\thm}[1]{Theorem \ref{#1}}
\newcommand{\prop}[1]{Proposition \ref{#1}}
\newcommand{\lma}[1]{Lemma \ref{#1}}

\newcommand{\cor}[1]{Corollary \ref{#1}}

\newcommand{\rarrow}{\rightarrow}

\newcommand{\fd}{\operatorname{fd}}

\newcommand{\chr}{\operatorname{char}}
\newcommand{\ho}{\operatorname{Ho}}
\newcommand{\spec}{\operatorname{Spec}}
\newcommand{\tor}{\operatorname{Tor}}

\newcommand{\calL}{{\mathcal L}}

\newcommand{\Alam}{{\mathcal Alg}_{\Lambda}}

\begin{document}

\title[On Simplicial Algebras with Noetherian Homotopy]
{On Simplicial Commutative Algebras with Noetherian Homotopy}
\author{James M. Turner}
\address{Department of Mathematics\\
Calvin College\\
3201 Burton Street, S.E.\\
Grand Rapids, MI 49546}
\email{jturner@calvin.edu}
\thanks{Research was partially supported by a grant from the National 
Science Foundation (USA)}
\date{\today}



\begin{abstract}
In this paper, we introduce a strategy for studying simplicial 
commutative algebras over general commutative rings $R$.  Given 
such a simplicial algebra $A$, this strategy involves replacing $A$ 
with a connected simplicial commutative $k(\wp)$-algebra $A(\wp)$, for 
each $\wp \in \spec(\pi_{0}A)$, which we call the {\bf connected 
component of $A$ at $\wp$}. These components retain most of the 
Andr\'e-Quillen homology of $A$ when the coefficients are 
$k(\wp)$-modules ($k(\wp)$ = residue field of $\wp$ in $\pi_{0}A$). 
Thus these components should carry quite a bit of the homotopy theoretic 
information for $A$. Our aim will be to apply this strategy to 
those simplicial algebras which possess {\bf Noetherian homotopy}. 
This allows us to have sophisticated techniques from commutative 
algebra at our disposal. One consequence of our efforts will be to 
resolve a more general form of a conjecture of Quillen 
that was posed in \cite{Tur}. 
\end{abstract}

\maketitle

\section*{Overview}
Our focus, in this paper, is to take the view that the study of 
Noetherian rings and algebras through homological methods is a special 
case of the study of simplicial commutative algebras having Noetherian 
homotopy type. Our goal is to show that such simplicial algebras can 
be given a suitably rigid structure in the homotopy category, which 
then allows us to bring in methods from commutative algebra. 
Such methods should enable more facile techniques from 
homological algebra to be ferried in for the purpose of elaborating 
the global structure of such simplicial algebras. 

To begin, we define for a simplicial 
commutative algebra $A$ to have {\it Noetherian homotopy} provided:
\begin{enumerate}
	
	\item $\pi_{0}A$ is a Noetherian ring, and
	
	\item each $\pi_{m}A$ is a finite $\pi_{0}A$-module.
	
\end{enumerate}

If, more strongly, $\pi_{*}A$ is a finite graded $\pi_{0}A$-module, we 
that $A$ has {\it finite Noetherian homotopy}.

In order to achieve a more systematic study of simplicial algebras with 
Noetherian homotopy, particularly to allow us a straighter path to 
proving our main result, Theorem B below, we first seek to rigidify 
the action of $\pi_{0}$ from the homotopy groups to the simplicial 
algebra. This is accomplished by the following:

\bigskip

\noindent {\bf Theorem A:} {\it \, Any simplicial commutative algebra 
$A$ is weakly equivalent to a connected simplicial supplemented 
$\pi_{0}A$-algebra.}

\bigskip

Theorem A provides the means to import in methods from commutative 
algebra, most notably localizations and completions. In particular, we 
use these methods as a means to provide a proof of a
conjecture posed in \cite{Tur} which generalizes a conjecture of 
Quillen regarding the vanishing of Andr\'e-Quillen homology. 
Our larger interests lie in providing an understanding of the 
of the homotopy type of a simplicial commutative 
algebra $A$ with Noetherian homotopy over a Noetherian ring $R$ 
through its Andr\'e-Quillen homology $D(A|R;-)$.  
Here we shall view this homology as a functor of $\pi_{0}A$-modules. 
This enables us to be specific about the homology's rigidity 
properties.

Before stating our result, we first need a homotopy invariant notion 
of complete intersection. To obtain one, we first define a map $A\to B$ of 
simplicial commutative $R$-algebras, augmented over a field $\ell$, 
to be {\it virtually acyclic} provided $D_{\geq 1}(B|A;\ell) = 0$. Also, 
if $W$ is a graded $\ell$-module, define the simplicial 
$\ell$-algebra $S_{\bullet}(W)$ by $$S_{\bullet}(W) = 
\bigotimes_{n}S(W_{n},n)$$ where $S(V,n)$ is the free commutative 
$\ell$-algebra generated by the Eilenberg-MacLane space $K(V,n).$

Define a simplicial commutative $R$-algebra $A$ over $\ell$ to be 
a {\it homotopy n-intersection,} for $n\geq 1$, provided there
is a commutative diagram
    $$
    \begin{array}{ccc}
    R & \stackrel{}{\longrightarrow} & R^{\prime} \\[1mm]
    \eta \downarrow \hspace*{10pt}
    &&
    \hspace*{10pt} \downarrow \eta^{\prime} \\[1mm]
    A & \stackrel{}{\longrightarrow} & A^{\prime}\\[1mm]
    \hspace*{0pt} \downarrow  \hspace*{0pt}
    &&
    \hspace*{0pt}  \downarrow  \\[1mm]
    \ell & \stackrel{=}{\longrightarrow} & \ell 
    \end{array}
    $$
with the horizontal maps being virtually acyclic over $\ell$ and 
in the homotopy category there is an isomorphism
$$A^{\prime}\otimes_{R^{\prime}}^{{\bf L}}\ell \cong 
S_{\bullet}(W)$$ with $W$ a graded $\ell$-module satisfying $W_{>n}=0.$
We call a general simplicial commutative $R$-algebra $A$ a {\it locally 
homotopy n-intersection} if, for 
each $\wp \in \spec (\pi_{0}A)$, $A$ is a
homotopy n-intersection over the residue field $k(\wp)$ 

Recall that the {\it flat dimension} 
of an $R$-module $M$ to be the positive integer $\fd_{R} M$ such that
\begin{equation}
	\fd_{R}M\leq m \Longleftrightarrow \tor^{R}_{i}(M,-) = 0 
	\quad \mbox{for} \quad i>m.
\end{equation}

\bigskip

\noindent {\bf Theorem B:}{\em \, Let $A$ be a simplicial 
commutative $R$-algebra with finite Noetherian 
homotopy, $\chr(\pi_{0}A)\neq 0$, and  
$\fd_{R} (\pi_{*}A)$ is finite. Then $D_{s}(A|R;-)=0$ for $s\gg 0$ 
if and only if $A$ is a locally homotopy 1-intersection.}

\bigskip

This resolves a conjecture posed in \cite{Tur} generalizing a 
conjecture of Quillen \cite[5.7]{Qui2}. 

\bigskip

\noindent {\bf Notes:}
\begin{enumerate}
	
	\item Theorem B fails when $\chr (\pi_{0}A) = 0$, as shown in 
	\cite{Tur}.
	
	\item Theorem B fails for general simplicial algebras having 
	Noetherian homotopy. The case of the simplicial algebras 
	$S(V,n)$ over a field of non-zero 
	characteristic provide counterexamples, by computations of Cartan 
	\cite{Car}.
	
	\item A homomorphism between Noetherian rings is a locally 
	complete intersection if and only if it is a locally homotopy 
    1-intersection, as shown in \cite{Avr, Tur}.
	
\end{enumerate}

Quillen further conjectured a more general result \cite[5.6]{Qui2} 
which drops the finite flat dimension condition. We would like to 
indicate a possible simplicial version of this conjecture of 
Quillen. To formulate it, we first indicate a special vanishing result
for Andr\'e-Quillen homology that we will prove.

\bigskip

\noindent {\bf Theorem C:} {\em \, Let $A$ be a simplicial 
commutative $R$-algebra with Noetherian homotopy. Then $D_{s}(A|R;-) = 
0$ for $s\geq 3$ if and only if $A$ is a locally homotopy 
2-intersection.}

\bigskip

This now leads us to pose the following:

\bigskip

\noindent {\bf Conjecture:}{\em \, Let $A$ have finite Noetherian 
homotopy with $\chr(\pi_{0}A)\neq 0$. Then $D_{s}(A|R;-)=0$ for $s\gg 0$ 
implies that $A$ is a locally homotopy 2-intersection.}
 
\bigskip

The strategy for proving Theorem B is to show that 
$D_{s}(A|R;k(\wp)) = 0$ for $s\geq 2$ for each $\wp \in \spec 
(\pi_{0}A)$. This is sufficient by a result of Andr\'e \cite[S.30]{And}. 
Following a strategy of Avramov \cite{Avr}, we use Theorem A coupled 
with commutative algebra techniques developed in \cite{AFH} to 
replace $A$ with $A(\wp)$, its {\it connected component at $\wp$},  
which has the following properties:
\begin{enumerate}
	
	\item $A(\wp)$ is a connected simplicial supplemented $k(\wp)$-algebra;
	
	\item $\fd_{R}(\pi_{*}A) <\infty$ implies that $A(\wp)$ has finite Noetherian 
	homotopy;
	
	\item $D_{s}(A|R;k(\wp)) \cong D_{s}(A(\wp)|k(\wp);k(\wp))$ for 
	$s\geq 2$.
	
\end{enumerate}
Theorem B now follows from the algebraic version of a theorem of Serre 
established in \cite{Tur}.

\subsection*{Acknowledgements} The author wishes to thank Lucho 
Avramov for sharing his expertise on commutative algebra and to Paul 
Goerss for sharing his expertise on Postnikov systems.

\setcounter{section}{0}

\section{Postnikov Systems and Theorem A}
Throughout this paper, we fix a commutative ring with unit $\Lambda$ and let 
$\Alam$ be the category of (unitary) commutative rings augmented 
over $\Lambda$. Finally, we denote by $_{\Lambda}\Alam$ the 
category of $\Lambda$-algebras in $\Alam$.

We will also be assuming the reader has an acquaintance with closed 
(simplicial) model category theory. Our main resource is \cite{Qui1}.
We will further need specific results on the model category structure 
for simplicial commutative rings and algebras. Our primary sources 
are \cite{Qui1,Qui3,Goerss}.

\subsection{Postnikov Systems}

Let $A$ be an object in the category $s\Alam$ of simplicial commutative 
rings over $\Lambda$. We 
review the construction of a Postnikov tower for $A$ derived from 
\cite{BDG, GT} which we will be use in the proof of Theorem A. 

Following \cite[\S 5]{GT}, define the {\it n}th Postnikov section of 
$A$ as follows: for fixed $k$, let $I_{n,k}\to A_{k}$ be the kernel of 
the map $$d:A_{k}\to \prod_{\phi:[m]\to[k]}A_{n}$$where $\phi$ runs 
over all injections in the ordinal number category with $m\leq n$, 
$d$ is induced by the maps $\phi^{*}: A_{k}\to A_{m}$, and $\prod$ 
denotes the product in the category of algebras augmented over 
$\Lambda$. Define
\begin{equation}\label{postdef}
	A(n)_{k} = A_{k}/I_{n,k}
\end{equation}	
Notice that there is a quotient map in $s\Alam$, $A \to A(n)$, and that 
if $k\leq n$, $A(n)_{k} = A_{k}$. There are also quotient maps
\begin{equation}\label{qmaps}
	q_{n} : A(n) \to A(n-1)
\end{equation}
and $A \cong \lim A(n)$. Let $F(n)$ be the fibre of $q_{n}$, i.e.
\begin{equation}\label{fibdef}
	F(n) = \ker (q_{n} : A(n) \to A(n-1)).
\end{equation}
Note that $F(n) \to A(n) \stackrel{q_{n}}{\to} A(n-1)$ forgets to 
a fibration sequence as simplicial abelian groups. As such, the 
following can be proved just as in \cite[5.5]{GT}.

\begin{lemma}\label{fibcomp}
	The homotopy groups of $F(n)$ are computed as follows:
    $$
    \pi_{k}F(n) =
    \begin{cases}
       \pi_{n}A & \quad k = n; \\
       0 & \quad k\neq n.
    \end{cases}
    $$
\end{lemma}

\subsection{Eilenberg-MacLane objects}

Following \cite[\S 5]{BDG}, define an object $A$ of $s\Alam$ to be 
{\it of type} $K_{\Lambda}$ if $\pi_{0}A \cong \Lambda$ and the 
higher homotopy groups of $A$ are trivial. Suppose $M$ is a 
$\Lambda$-module. We say that a map $A\to B$ is {\it of type} 
$K_{\Lambda}(M,n)$ $n\geq 1$, if $A$ is of type $K_{\Lambda}$, 
$\pi_{0}B \cong \Lambda$, $\pi_{n}B \cong M$ (as a $\Lambda$-module), 
all other homotopy groups of $B$ are trivial, and the map $A\to B$ is 
a $\pi_{0}$-isomorphism. 

For a general map $f : A \to B$ in $s\Alam$, let $C$ be the pushout of 
the diagram $B^{\prime} \leftarrow A^{\prime} \rightarrow 
A(0)^{\prime}$ obtained by using a functorial construction to 
replace $A$ by a cofibrant object and the two maps $A \to B$ and $A 
\to A(0)$ by cofibrations. There is then a commutative diagram
\begin{equation}\label{postdiag}
\begin{array}{ccc}
A & \stackrel{f}{\longrightarrow} & B \\[1mm]
\sim \uparrow \hspace*{10pt}
&&
\hspace*{10pt} \uparrow \sim \\[1mm]
A^{\prime} & \stackrel{f^{\prime}}{\longrightarrow} & B^{\prime}\\[1mm]
\hspace*{0pt}  \downarrow  \hspace*{0pt}
&&
\hspace*{0pt}  \downarrow \\[1mm]
A(0)^{\prime} & \stackrel{\Delta_{n}(f)}{\longrightarrow}
& C(n+1)
\end{array}
\end{equation}
The bottom map $\Delta_{n}(f)$ is called the {\it difference 
construction of f.} The following can be proved just as in  
\cite[6.3]{BDG}.

\begin{proposition}\label{kinv}
	Suppose that $A \to B$ is a map of simplicial commutative algebras 
	which is a $\pi_{0}$-isomorphism and whose homotopy fibre $F$ is 
	(n-1)-connected. Let $M = \pi_{n}F$. Then $M$ is naturally a 
	$\Lambda$-module for $\Lambda = \pi_{0}B$ and $\Delta_{n}(f)$ is a 
	map of type $K_{\Lambda}(M,n+1)$. If $\pi_{k}F$ vanishes except for 
	$k = n$, then the right-hand square in \ref{postdiag} is a homotopy 
	fibre square.
\end{proposition}	

\subsection{Differentials functor}

For an object $A$ in $\Alam$, define its {\it 
$\Lambda$-differentials} 
to be the $\Lambda$-module $$D_{\Lambda}A = J/J^{2}\otimes_{A}\Lambda$$
where $J$ is the kernel of the product $A\otimes A \to A$. 
As a functor to the category of 
$\Lambda$-modules, $D_{\Lambda}$ posseses a right adjoint - the 
functor $$(-)_{+} : Mod_{\Lambda} \to \Alam$$ defined by $M_{+} = 
M\oplus\Lambda$ with the usual twisted product
$$
(x, a)\cdot (y, b) = (bx+ay, ab).
$$
An equivalent identification of the differentials functor
\begin{equation}\label{diffind}
	D_{\Lambda} \cong I/I^{2}\otimes_{A}\Lambda,
\end{equation}
where $I$ is the augmentation ideal of $A$, which can be 
seen to follow from Yoneda's lemma. 

The next proposition is proved in \cite[\S II.5]{Qui1}.

\begin{proposition}\label{homprop}
	The prolonged adjoint pair of functors
	$$
	D_{\Lambda} : s\Alam \Longleftrightarrow sMod_{\Lambda} : (-)_{+}
	$$
	induces an adjoint pair on the homotopy categories
	$$
	{\bf L}D_{\Lambda} : \mbox{Ho}(s\Alam) \Longleftrightarrow 
	\mbox{Ho}(sMod_{\Lambda}) : {\bf R}(-)_{+}.
	$$
\end{proposition}

Finally, the following useful property of the derived functor of differentials 
follows from \cite[7.3]{Qui3}.

\begin{proposition}\label{indprop}
	If $f: A \to B$ is a $\pi_{\leq n}$-isomorphism, 
	then ${\bf L}D_{\Lambda}(f)$ is a $\pi_{\leq n}$-isomorphism.
\end{proposition}

\subsection{Characterizing $K_{\Lambda}(M,n)$-type}

Fix a $\Lambda$-module $M$. In $sMod_{\Lambda}$, the fibration 
$p_{n}: E(M,n) \to K(M,n)$ is determined by the Dold-Kan correspondence by 
to correspond to the map of normalized chain complexes 
$\{M\stackrel{1}{\to} M\} \to \{M\}$ with the source concentrated in 
degrees $n$ and $n-1$, the target concentrated in degree $n$, and the 
map being the identity in degree $n$ and trivial otherwise.

Applying $(-)_{+}$ to $p_{n}$ gives a $K_{\Lambda}(M,n)$-type fibration 
in $s\Alam$
$$
(p_{n})_{+} : E_{\Lambda}(M,n) \rarrow K_{\Lambda}(M,n)
$$
which we call the {\it canonical map of type $K_{\Lambda}(M,n)$}. 

\begin{proposition}\label{klamprop}
	Let $A \to B$ be of type $K_{\Lambda}(M,n)$ between cofibrant objects 
	in $s\Alam$. Then there is a commuting diagram in $s\Alam$
	$$
    \begin{array}{ccc}
    A & \stackrel{\sim}{\longrightarrow} & E_{\Lambda}(M,n)\\[1mm]
    \hspace*{10pt}  \downarrow  \hspace*{0pt}
    &&
    \hspace*{0pt}  \downarrow p_{n}\\[1mm]
    B & \stackrel{\sim}{\longrightarrow} & K_{\Lambda}(M,n) 
    \end{array}
	$$
	with the horizontal maps being weak equivalences.
\end{proposition}	

\noindent {\it Proof.} To begin, note that the canonical map $B \to \Lambda$ 
is (n-1)-connected. Thus the induced map $D_{\Lambda}B \to 0$ is 
(n-1)-connected by \prop{indprop}. Let $I = \ker (B \to \Lambda)$.  
Filtering $B$ by powers of $I$ we note that $B$ cofibrant implies that
$$
I^{q}/I^{q+1} = S^{\Lambda}_{q}(I/I^{2}) \cong 
S^{\Lambda}_{q}(D_{\Lambda}B)
$$
where the last identity always holds when the augmentation is 
surjective, by (\ref{diffind}). Thus there is a convergent 
spectral sequence
$$
E^{1}_{p,q} = H_{p+q}[S^{\Lambda}_{q}(D_{\Lambda}B)] \Longrightarrow 
\pi_{p+q}B.
$$
From the connectivity indicated above and \cite[7.40]{Qui3}, 
$E^{1}_{p,q} = 0$ for $0<p+q\leq 2(q-2)+n$. Thus we obtain
$$
M \cong \pi_{n}B \cong \pi_{n}D_{\Lambda}B.
$$
Thus there is an n-connected map $D_{\Lambda}B \to K(M,n)$ and its 
adjoint $B \to K_{\Lambda}(M,n)$ will be a weak equivalence by the
computations above and the assumption that $A \to B$ is of type 
$K_{\Lambda}(M,n)$. 

Finally, $A \to \Lambda$ is a weak equivalence, hence $D_{\Lambda}A 
\to 0$ is a weak equivalence by \prop{indprop}. Since $A$, and hence 
$D_{\Lambda}A$, are cofibrant, the composite $D_{\Lambda}A \to 
D_{\Lambda}B \to K(M,n)$ lifts to a map $D_{\Lambda}A \to 
E(M,n),$ whose adjoint $A \to E_{\Lambda}(M,n)$ is necessarily a weak 
equivalence. \hfill $\Box$  

\subsection{Proof of Theorem A}

Fix an object $A$ in $s\Alam$. We will show, by induction, that 
there is a map $X \to Y$ in $s _{\Lambda}\Alam$ and a commutative 
diagram in $\mbox{Ho}(s\Alam)$
\begin{equation}\label{indstep}
    \begin{array}{ccc}
    A(n) & \stackrel{\sim}{\longrightarrow} & X\\[1mm]
    \hspace*{0pt} q_{n} \downarrow  \hspace*{10pt}
    &&
    \hspace*{0pt}  \downarrow \\[1mm]
    A(n-1) & \stackrel{\sim}{\longrightarrow} & Y
    \end{array}	
\end{equation}	
with the horizontal maps being equivalences. It is clear for $n = 0$ 
as $A(0) \to \Lambda$ is a weak equivalence. 

Using \ref{postdiag}, some closed model category theory and induction, we may 
assume that 
there is a trivial fibration $\sigma : A(n-1)^{\prime} \to Y$ with the target $Y$ a 
cofibrant object in $s _{\Lambda}\Alam$. 

\begin{lemma}\label{speccase}
	Let $M = \pi_{n}A$. Then there is a commuting diagram in 
	$\mbox{Ho}(s\Alam)$ of the form
    $$
	\begin{array}{ccc}
    A(n-1)^{\prime} & \stackrel{}{\longrightarrow} & C(n+1)\\[1mm]
    \hspace*{10pt} \sim \downarrow \sigma \hspace*{15pt}
    &&
    \hspace*{0pt}  \downarrow \sim \\[1mm]
    Y & \stackrel{}{\longrightarrow} & K_{\Lambda}(M,n+1)
    \end{array}	
	$$
	with the top arrow from \ref{postdiag}.
\end{lemma}	

\noindent {\it Proof.} First, note that since $\sigma : A(n-1)^{\prime} \to Y$ is a 
trivial fibration between suitably cofibrant objects (see above) it 
follows from that and from \ref{diffind} that
$$
D_{\Lambda}\sigma : D_{\Lambda}A(n-1)^{\prime} \to D_{\Lambda}Y
$$
is a trivial fibration between cofibrant objects in $sMod_{\Lambda}$. By 
\cite[I.1.7]{Qui1}, $D_{\Lambda}\sigma$ has a homotopy left inverse $i$ 
($i\circ D_{\Lambda}\sigma \simeq \mbox{Id}_{D_{\Lambda}A(n-1)}$). 

Next, utilizing \lma{klamprop}, let 
$t : A(n-1)^{\prime} \to K_{\Lambda}(M,n+1)$ be the 
composite of $A(n-1)^{\prime} \to C(n+1) \to K_{\Lambda}(M,n+1)$. 
Let $w : D_{\Lambda}Y \to K(M,n+1)$ be the composite 
$(D_{\Lambda}t)\circ i$. Then $w\circ D_{\Lambda}\sigma \simeq 
D_{\Lambda}t$ and the result now follows from \prop{homprop}.
\hfill $\Box$

\bigskip

From the previous lemma, we may form the homotopy pullback diagram in $s 
_{\Lambda}\Alam$
\begin{equation}\label{diag1}
	\begin{array}{ccc}
    X & \stackrel{}{\longrightarrow} & E_{\Lambda}(M,n+1)\\[1mm]
    \hspace*{10pt} \downarrow  \hspace*{10pt}
    &&
    \hspace*{20pt}  \downarrow (p_{n})_{+} \\[1mm]
    Y & \stackrel{}{\longrightarrow} & K_{\Lambda}(M,n+1).
    \end{array}	
\end{equation}
By \prop{kinv}, the diagram below is also a homotopy pullback in 
$s\Alam$  
\begin{equation}\label{diag2}
	\begin{array}{ccc}
    A(n)^{\prime} & \stackrel{}{\longrightarrow} & A(0)^{\prime}\\[1mm]
    \hspace*{0pt} q_{n}^{\prime} \downarrow  \hspace*{20pt}
    &&
    \hspace*{20pt}  \downarrow \Delta[q_{n}]  \\[1mm]
    A(n-1)^{\prime} & \stackrel{}{\longrightarrow} & C(n+1).
    \end{array}	
\end{equation}
By \prop{klamprop} and \lma{speccase}, there is an induced map of 
diagrams \ref{diag2} to \ref{diag1} in the category $\mbox{Ho}(s\Alam)$.
Since fibrations and pullbacks in $s\Alam$ are fibrations and pullbacks 
as simplicial groups, a computation of homotopy groups can be 
performed utilizing \lma{fibcomp} to show that the induced map 
$A(n)^{\prime} \to X$ is a weak equivalence. This completes the induction step. 

\section{Andr\'e-Quillen homology and Theorems B and C}

\subsection{Base change property of Andr\'e-Quillen homology}

Recall that the {\it cotangent complex} of a simplicial $R$-algebra 
$A$ is defined to be the object of $\ho (Mod_{A})$
\begin{equation}\label{ccdef}
\calL(A|R) := \Omega_{P|R}\otimes_{P}A
\end{equation}
where the $T$-module $\Omega_{T|S} = J/J^{2}$, $J = \ker (T\otimes_{S}T \to T)$, 
denotes the {\it Kahler differentials} of an $S$-algebra $T$, and
$P \to A$ is a cofibrant replacement of $A$ as a simplicial 
$R$-algebra. 

\bigskip

\noindent {\bf Note:} As in \S 1.3, $\Omega_{T|S}$ is left adjoint to 
the functor $M \mapsto M\oplus T$ where the image has a $T$-algebra 
structure with $M^{2} = 0$.

\bigskip 

Also recall that given another simplicial $R$-algebra 
$B$, the {\it derived tensor product} of $A$ and $B$ to be the object 
of $\ho (sMod_{R})$
$$
A\otimes_{R}^{{\bf L}}B := P\otimes_{R}Q
$$
where $Q \to B$ is a cofibrant replacement of $B$. 

We now derive a base change property for the cotangent complex 
following \cite{Qui3}.

\begin{lemma}\label{l1}
	If $\tor^{R}_{q}(A_{k},B_{k}) = 0$ for all $k\geq 0$ and all $q>0$ 
	then $A\otimes^{{\bf L}}_{R}B \simeq A\otimes_{R} B$.
\end{lemma}

\noindent {\it Proof.} This follows immediately from the spectral 
sequence \cite[\S II.6]{Qui1}
$$
E^{2}_{p,q} = \pi_{p}\tor_{q}^{R}(A,B) \Longrightarrow 
\pi_{p+q}(A\otimes_{R}^{{\bf L}}B).
$$
\hfill $\Box$

\bigskip

\begin{lemma}\label{l2}
	$\Omega_{A\otimes_{R}B|B}\cong \Omega_{A|R}\otimes_{R}B$
\end{lemma}

\noindent {\it Proof.} Let $A^{\prime} = A\otimes_{R}B$ and fix an 
$A^{\prime}$-module $M$. Then
$$
\begin{array}{cl}
	\hom_{A^{\prime}}(\Omega_{A^{\prime}|B},M) & \cong 
	\hom_{_{B}Alg_{A^{\prime}}} (A^{\prime}, M\oplus A^{\prime})\\ 
	 & \cong \hom_{_{R}Alg_{A}} (A, M\oplus A) \\
	 & \cong \hom_{A} (\Omega_{A|R}, M) \\
	 & \cong \hom_{A^{\prime}} (\Omega_{A|R}\otimes_{R}B,M).\\
\end{array}
$$
The result now follows from Yoneda's lemma. \hfill $\Box$
	
\bigskip 

\begin{proposition}\label{bchange}
	$\calL(A\otimes_{R}^{{\bf L}}B|B) \simeq \calL(A|R)\otimes^{{\bf 
	L}}_{R}B$
\end{proposition}

\noindent {\it Proof.} Fix cofibrant replacements $P$ and $Q$ for 
$A$ and $B$, respectively. Then
\begin{equation}\label{r1}
	\calL(A\otimes^{{\bf L}}_{R}B|B) = \Omega_{P\otimes_{R}Q|Q} \cong 
	\Omega_{P|R}\otimes_{R}Q
\end{equation}
by \lma{l2}. Since $P$ is projective as a simplicial $R$-module then 
$\Omega_{P|R}$ is a projective $P$-module. Thus, by \lma{l1}, the map 
$\Omega_{P|R} \stackrel{\sim}{\to} \Omega_{P|R}\otimes_{P}A$ is a weak equivalence.
Since $Q$ is projective, \lma{l1} further tells us that
\begin{equation}\label{r2}
	\Omega_{P|R}\otimes_{R}Q \stackrel{\sim}{\rarrow} 
	(\Omega_{P|R}\otimes_{P}A)\otimes_{R}Q \cong \calL(A|R)\otimes^{{\bf 
	L}}_{R}B
\end{equation}
is a weak equivalence. The result now follows by combining \ref{r1} 
with \ref{r2}. \hfill $\Box$ 

\bigskip

\begin{corollary}\label{bchange2}
	As a functor of $A\otimes_{R}B$-modules, $D_{*}(A\otimes^{{\bf 
	L}}_{R}B|B;-) \cong D_{*}(A|R;-).$
\end{corollary}

\noindent {\it Proof.} This follows from \prop{bchange} and the 
identity $D_{*}(T|S;M) := \\ \pi_{*}[\calL(T|S)\otimes_{T}M]$. 
\hfill $\Box$

\bigskip

\subsection{Proof of Theorem B}

We first recall the main result of \cite{Tur}.

\begin{theorem}\label{algserre}
	Let $A$ be a homotopy connected simplicial supplemented commutative algebra 
	over a field $\ell$ of non-zero characteristic. Then 
	$D_{s}(A|\ell;\ell) = 0$ for $s\gg 0$ implies that there is an 
	equivalence $S_{\ell}(D_{1}(A|\ell;\ell),1) \cong A$ in the homotopy 
	category.
\end{theorem}	

We now begin by establishing a special case of Theorem A. To that end
let $A$ be a simplicial commutative $R$-algebra and assume 
that the unit $R \to \pi_{0}A = \Lambda$ is a surjection. For $\wp \in \spec 
\Lambda$, define the {\it connected component of $A$ at $\wp$} to be 
the connected simplicial supplemented $k(\wp)$-algebra
$$
A(\wp) = A\otimes^{{\bf L}}_{R}k(\wp).
$$

\begin{lemma}\label{ccatwp}
	Let $A$ be as above. Then
\begin{enumerate}
		
	\item $D_{*}(A|R;k(\wp)) \cong D_{*}(A(\wp)|k(\wp);k(\wp))$, and

	\item if $A$ also has finite Noetherian homotopy and 
	$\fd_{R}(\pi_{*}A)<\infty$ it follows that $A(\wp)$ has finite Noetherian 
	homotopy.
	
\end{enumerate}

\end{lemma}

\noindent {\it Proof.} 1. follows from \cor{bchange2}. For 2., 
\cite[\S II.6]{Qui1} gives a spectral sequence
$$
E^{2}_{s,t}=\tor^{R}_{s}(\pi_{t}A, k(\wp)) \Longrightarrow 
\pi_{s+t}(A\otimes^{{\bf L}}_{R}k(\wp)).
$$
From the finiteness conditions, each $E^{2}_{s,t}$ is a finite 
$k(\wp)$-module and vanishes for $s,t \gg 0$. Thus 
$A\otimes^{{\bf L}}_{R}k(\wp))$ has finite Noetherian homotopy. 
\hfill $\Box$ 

\begin{corollary}\label{cor}
	Let $A$ be as in \lma{ccatwp}.2 and further assume that $\chr (k(\wp)) 
	\neq 0$. Then $D_{s}(A|R;k(\wp))=0$ for $s\gg 0$ 
    implies that $D_{s}(A|R;k(\wp))=0$ for $s\geq 2$.
\end{corollary}

\noindent {\it Proof.} This follows from \lma{ccatwp} and 
\thm{algserre}. \hfill $\Box$

\bigskip

Now assume that the simplicial algebra $A$ in question 
is a homotopy connected simplicial supplemented $\Lambda$-algebra, by Theorem A.
We further assume that $A$ has Noetherian homotopy. 

Fix $\wp \in \spec \Lambda$ and let $\widehat{(-)}$ denote the 
completion functor on $R$-modules at $\wp$. Define the homotopy 
connected simplicial supplemented $\widehat \Lambda$-algebra $A^{\prime}$ by
$$
A^{\prime} = A\otimes^{{\bf L}}_{\Lambda}\widehat \Lambda.
$$

\begin{proposition}\label{final}
	Suppose $A$ is a simplicial commutative $R$-algebra, with $R$ a Noetherian ring. 
	Then $\pi_{*}A^{\prime} \cong \widehat{\pi_{*}A}$ and
	there exists a (complete) Noetherian $R^{\prime}$ that fits into the 
	following commutative diagram in $\ho (s_{R}{\mathcal Alg})$
$$
\begin{array}{ccc}
R & \stackrel{\eta}{\longrightarrow} &
A \\[1mm]
\phi \downarrow \hspace*{10pt}
&&
\hspace*{10pt} \downarrow \psi \\[1mm]
R^{\prime} & \stackrel{\eta^{\prime}}{\longrightarrow} & A^{\prime}
\end{array}
$$
     with the following properties:
	 \begin{enumerate}
		 
		 \item $\phi$ is a flat map and its closed fibre $R^{\prime}/\wp 
		 R^{\prime}$ is weakly regular;
		 
		 \item $\psi$ is a $D_{*}(-|R;k(\wp))$-isomomorphism;
		 
		 \item $\eta^{\prime}$ induces a surjection $\eta^{\prime}_{*}: R^{\prime} \to 
		 \pi_{0}A^{\prime}$;
		 
		 \item $\fd_{R} (\pi_{*} A)$ finite implies that $\fd_{R^{\prime}} 
		 (\pi_{*} A^{\prime})$ is finite
		 
	\end{enumerate}

\end{proposition}	

\noindent {\it Proof:} First, Quillen's spectral sequence 
\cite[II.6]{Qui1} $\tor_{*}^{\Lambda}(\pi_{*}A,\widehat \Lambda) 
\Longrightarrow \pi_{*}A^{\prime}$ collapses to give the first result 
since $\widehat \Lambda$ is flat over $\Lambda$ and each $\pi_{m}A$ 
is finite over $\Lambda$ \cite[8.7 and 8.8]{Mat}.

Next, by \cite[1.1]{AFH}, the unit ring homomorphism $R 
\to \widehat \Lambda$ factors as $R \stackrel{\phi}{\rarrow} R^{\prime} 
\stackrel{\eta^{\prime}_{*}}{\rarrow} \widehat \Lambda$ with $\phi$ 
having the properties described in 1. and $\eta^{\prime}_{*}$ is a 
surjection. Thus the induced map $\eta^{\prime}: R^{\prime} \to 
A^{\prime}$ induces a surjection on $\pi_{0}$, giving 3., and the 
desired diagram commutes. 

Now, by the transitivity sequence \cite[4.12]{Qui3} applied to $R \to 
A \to A^{\prime}$, 2. follows from the isomorphism
$$
D_{*}(A^{\prime}|A;k(\wp)) \cong D_{*}(\widehat \Lambda 
|\Lambda;k(\wp)) \cong 0
$$
which follows from \cor{bchange2}.

Finally, 4. follows from \cite[3.2]{AFH}, as $A$ has Noetherian 
homotopy.  \hfill $\Box$ 

\bigskip 

Now, let $A$ have finite Noetherian homotopy with $D_{s}(A|R;-) = 0$ 
for $s\gg 0$. From \prop{final}, \thm{algserre}, \cor{cor}, and 
\cite[\S S.30]{And}, if $\fd_{R} 
(\pi_{*}A)< \infty$ then $A(\wp) \cong S_{k(\wp)}(D_{1}(A|R;k(\wp),1)$, for each 
$\wp\in\spec (\pi_{0}A)$, if and only if $D(A|R;-) = 0$. Thus Theorem B follows 
from the definition of locally homotopy complete intersection (see 
introduction) and a transitivity sequence argument.

\subsection{Proof of Theorem C}

Let $A$ be a simplicial commutative $R$-algebra with Noetherian homotopy. 
It follows from \lma{ccatwp}.1, \prop{final}, and \cite[\S S.30]{And}, 
that $D_{\geq 3}(A|R;-) = 0$
if and only if $D_{\geq 3}(A(\wp)|k(\wp);k(\wp)) = 0$, for all 
$\wp\in\spec (\pi_{0}A)$. From the definition of locally virtual 
homotopy complete intersection (see introduction), Theorem C will follow if we can show 
that, for each prime ideal $\wp$, $A(\wp) \cong S_{\bullet}(D_{\leq 
2}(A|R;k(\wp)))$ in the homotopy category. But this in turn follows 
from \cite[(2.2)]{Tur}.

\end{document}